\documentclass[12pt, oneside]{article}   	
\usepackage[letterpaper,margin=1in]{geometry}
\usepackage{graphicx}
\usepackage{amssymb}
\usepackage{amsmath}
\usepackage{color}
\usepackage{authblk}
\newcommand{\eps}{\varepsilon}
\newcommand{\R}{{\mathbb R}}
\newcommand{\Z}{{\mathbb Z}}
\newcommand\eJG[1]{{\color{black}{#1}}}

\title{The Legacy of the Cartwright-Littlewood Collaboration
}
\author{John Guckenheimer}
\affil{Mathematics Department \\ Cornell University \\ Ithaca, New York 14853}
\begin{document}

\maketitle

\begin{abstract}
Mary L. Cartwright and John E. Littlewood published a short ``preliminary survey" in 1945 describing results of their investigation of the \emph{forced van der Pol equation}
\begin{equation*}
\ddot{y}-k(1-y^2)\dot{y}+y = b \lambda k \cos(\lambda t+a)
\end{equation*}
in which $b,\lambda,k,a$ are parameters with $k$ large. 
\eJG{Their description of dynamical behavior now known as \emph{chaos} in this dissipative dynamical system was a landmark in dynamical systems theory.} 
Littlewood's ``monster paper" containing the details of their investigation finally appeared twelve years later in the journal Acta Mathematica. I review here the context in which Cartwright and Littlewood worked when they wrote their 1945 paper and the enduring mathematical legacy of their discoveries. I also give brief pointers to research they inspired in other application areas.
\end{abstract}

\section{Introduction}

Mary L. Cartwright and John E. Littlewood investigated the  forced van der Pol equation
\begin{equation}
\ddot{y}-k(1-y^2)\dot{y}+y = b \lambda k \cos(\lambda t+a)
\label{fvdp}
\end{equation}
following a 1938 memorandum from the Department of Scientific and Industrial Research Radio Research Board
requesting help from mathematicians in solving nonlinear differential equations~\eJG{\cite{MR160988,CL45}}. The equations of interest were second order nonlinear differential equation models of radios containing nonlinear thermionic valves (also known as vacuum tubes) as well as linear capacitors, inductances and resistors. Equation \eqref{fvdp} was one of several differential equations that had been studied extensively by \eJG{Balthazar van der Pol}  in the 1920's and 1930's. After van der Pol helped develop electronic devices as a physics student at Cambridge University, he became a leader in establishing international standards for radio while maintaining a strong interest in the mathematics of nonlinear ordinary differential equations \eJG{and their applications in diverse fields~\cite{israel2004technological}}. Cartwright and Littlewood cited two of van der Pol's papers as inspiration their work: a 1927 Nature letter \cite{VV27} by van der Pol and van der Mark which reports subharmonic outputs of different frequencies in an electronic circuit modeled by equation \eqref{fvdp}, and van der Pol's 1934 IRE paper \cite{Vpl34} on oscillatory solutions of nonlinear differential equations. 

Relying on earlier work of Birkhoff~\eJG{\cite{birkhoff1932quelques}}, Cartwright and Littlewood recognized that if equation \eqref{fvdp} had subharmonic solutions of different periods similar to those observed experimentally by van der Pol and van der Mark, the phase portrait of this dynamical system must be very complicated topologically. The term \emph{chaos}\footnote{\eJG{\emph{Chaos} has been given multiple mathematical definitions. The most stringent, which fits the phenomena described by Cartwright and Littlewood, is the presence of \emph{hyperbolic invariant sets}.  Informally these are often called \emph{horseshoes}, following the work of Smale~\cite{Smale65}. An example is displayed below in Figure 2.}}  is now used to describe solutions like those found by Cartwright and Littlewood. Chaos explains how a clockwork universe that follows trajectories of deterministic ordinary differential equations can still be unpredictable. The subject has attracted a colorful list of investigators from varied disciplines -- as recounted by James Gleick in his 1987 book ``Chaos: Making a New Science"~\cite{Gleick87} that has sold over a million copies. While Cartwright and Littlewood began their investigations in response to the Radio Research Board memorandum that sought closer ties between mathematicians and engineers, they worked within the boundaries of mathematics as a self contained subject. McMurran and Tattersall~\cite{MR1427114} described how the idiosyncratic rules they followed in their collaboration encouraged their transformative research. 

Papers without proofs are seldom published in mathematics research journals like the Journal of the London Mathematical Society, but World War II was an exceptional time, and chaos was already a phenomenon of extraordinary interest -- long before it acquired a mathematical definition. Norman Levinson noted in his Mathematical Reviews entry for the Cartwright-Littlewood paper that he had asked in \cite{MR11505} whether second order ordinary differential equations with periodic dependence on time could have complicated solutions like the ones found by Cartwright and Littlewood. He pursued their work by replacing the quadratic nonlinear term in the forced van der Pol equation with a discontinuous step function that makes the system piecewise linear. Levinson solved the modified system by concatenating explicit solutions of the linear systems and found parameter regimes with chaotic solutions~\cite{MR30079}. He showed further that the piecewise linear system could be perturbed to a smooth (high degree) polynomial system that retained chaotic solutions. These results established a rigorous foundation for the conceptual insights of Cartwright and Littlewood, though in an altered setting. Inspired by Levinson's work, Stephen Smale~\cite{Smale65}, discovered the \emph{horseshoe}, a still simpler discrete time dynamical system with chaotic trajectories. Smale then went much further, formulating an audacious research program in an article published in the Bulletin of the American Mathematical Society~\cite{Smale67} that led to to a geometric characterization of \emph{structurally stable}\footnote{A diffeomorphism $f$ is $C^r$-\emph{structurally stable} if $C^r$ perturbations of $f$ are all \emph{topologically equivalent}; i.e. there are homeomorphisms mapping trajectories  of perturbations to those of $f$.} smooth diffeomorphisms on compact manifolds~\cite{MR287580,MR932138}.

Fitting the forced van der Pol equation into the framework pioneered by Smale was still a challenge since it cannot be solved with explicit analytic formulas. While Littlewood's complicated 1957 papers~\cite{Littlewood57,Littlewood57b} gave further information about its solutions, they were far from a comprehensive description of the phase portrait of the equation at a single set of parameters or the bifurcations that occur in this phase portrait as system parameters are varied. The work described here filled that gap, using the following vector field on the manifold $\R^2 \times S^1$ that is equivalent to the equation \eqref{fvdp}~\cite{GHW03}
\begin{equation}
\begin{array}{rcl}
   \eps \dot{x} & = & y + x - \frac{x^3}{3}  \\
   \dot{y} & = & - x + a \sin(2\pi \theta) \\
   \dot{\theta} & = & \omega
\end{array}
\label{fvdpv}
\end{equation}
The symbols $a$ and $y$ have different meanings in equations \eqref{fvdp} and \eqref{fvdpv}. In the remainder of this paper, we refer to the vector field \eqref{fvdpv} and its flow with the acronym FVDP. The map $\phi$ that advances solutions for time  $2\pi/\omega$ maps cross-sections of constant $\theta$ to themselves, enabling the system to be reduced to a diffeomorphism of a cross-section of constant $\theta$, called a \emph{Poincar\'e map} or \emph{return map}. 

As described in the next section of this paper, the Radio Board memorandum sought efficient methods for finding stable periodic orbits and their periods within a small class of differential equations that includes FVDP. Determining which regions of the $(a,\omega,\eps)$ parameter space of FVDP have the desired oscillatory solutions is part of a \emph{bifurcation analysis} that partitions the parameter space into sets of topologically equivalent flows and gives a geometric description of the dynamics in each of these. Bifurcation analysis is an important technology for the design of systems that manifest desired dynamical behaviors like the stable periodic orbits of radio transmitters. FVDP became a prominent case study for developing methods for bifurcation analysis of dynamical systems with multiple time scales systems.

Sections 3,4 and 5 of this paper describe developments in three areas of dynamical systems theory that have been used to produce a more complete description of the dynamics of FVDP than was achieved by Cartwright and Littlewood:
\begin{itemize}
\item
Bifurcation theory of generic families of diffeomorphisms and flows,
\item
Geometric singular perturbation theory, and
\item
Numerical methods for continuation of invariant manifolds.
\end{itemize}
Drawing upon all three of these areas, the Nonlinearity paper of Haiduc~\cite{Haiduc09} proves that FVDP has parameter regions in which the system displays chaotic dynamics and is also structurally stable. Section 6 of this paper gives pointers to a few application areas that have benefitted from and contributed to these developments in nonlinear dynamics. 

\section{The Radio Board Memorandum: Oscillations}

\eJG{In her Presidential address to the London Mathematical Society~\cite{MR160988}, Cartwright refers to the ``memorandum from the Department of Scientific and Industrial Research dated 11 January, 1938 which was sent to the London Mathematical Society." I requested help from Caroline Series as editor of this special issue in finding a copy of this memorandum. She then sought assistance from June Barrow-Green who located a draft of the memorandum by F. Morley Colebrook~\cite{smith1954mr} entitled ``A note on certain types of Non-Linear Differential Equations involved in Radio Engineering" in the minutes of a Radio Research Board Committee stored in an off-site facility of the British National Archives. Those minutes have hand-written corrections in response to a letter from Sydney Chapman suggesting that the report be sent to university mathematics departments and members of the London Mathematical Society rather than the LMS as an organization. Cartwright also refers to her correspondence with Colebrook, van der Pol and Appleton that helped clarify the mathematical issues of most interest and importance for radio engineering. In a separate item from the draft memorandum, the Committee minutes state that their chief concern was the ``development of single frequency oscillators."}

The following three quotes from the draft memorandum state what the Radio Research Board Committee hoped mathematicians could do to help their efforts:
\begin{description}
\item
``The radio-engineer wishes to be able to find explicit analytical solutions for such equations, solutions having some of the comparative simplicity and utility as those available for linear differential equations he is already familiar with; that is to say, an explicit formulation of the dependent variable in terms of the independent variable and the parameters of the system, preferably a formulation in terms of known or tabulated functions in order that the solution and its variation with the different parameters may be evaluated in actual numbers without a prohibitive labour of computation.''
\item
``It is possible that one reason for the apparent lack of published work in this field on the part of mathematicians may be a deficiency of liaison between mathematicians on the one hand and engineers and physicists on the other. One of the objects of this note is at least to do something to meet that deficiency by calling attention to the fact that there is an actual and increasing need for really expert guidance in this matter. Another possible, indeed probable, reason is that solutions of such equations, in the sense indicated above, do not in fact exist. That is to say, it may well be that any explicit solutions, even if obtainable, would be of so complicated a form that they would be no more useful, as a source of exact information about the behaviour of the system, than the original differential equations themselves. Even in this case, however, the guidance and advice of the mathematician would be valuable, in preventing the waste of time and energy spent in pursuit of a will-o'-the wisp.''
\item
``In the majority of analytical problems involved in radio technique, the dependent variable is, in its final steady state, a periodic function of time. In general, considerably less interest attaches to the precise numerical evaluation of the solution than to a knowledge of the dependence of the various elements of the solution, particularly the fundamental frequency, on the various parameters. The magnitudes of these parameters are at the disposal of the engineer and experimenter, and for him the important objective of the analysis is information as to the appropriate values to give to them in relation to some practical objective in radio communication. This means that analytical solutions, even complicated and implicit solutions, are likely to be of greater practical value than any form of graphical or numerical solution in which the individual parameters have lost their identity, and in which the significance of individual parameters can only be arrived at by evaluating numerical solutions over wide ranges of assigned values.''
\end{description}

The fears expressed by the committee in these quotations were well founded because explicit formulas for solving FVDP do not exist. The problem highlighted in the memorandum of determining formulas for the period of oscillators (stable periodic orbits) as a function of system parameters has been studied for the (unforced) van der Pol equation with only partial success. The limit value of the periods of the van der Pol oscillations as its time scale parameter $\eps \to 0$ is readily determined with an explicit formula for the limit equation, but the $\eps $ dependence of the period is singular~\cite{Grasman87}. Different singular asymptotic expansions have been proposed, but it is likely that they all diverge and fail to give precise values for the periods. Furthermore, numerical estimates for the periods based on simulation are only as good as the accuracy of numerical integration algorithms, and these encounter increasing difficulty as $\eps$ approaches $0$. While techniques that utilize Laplace transforms and their inverses (called Borel transforms) have been successful in solving some singular perturbation problems in dynamical systems, the ones involving periods of the van der Pol relaxation oscillations remain open. Today, when computers are rapidly expanding the boundaries of machine intelligence with formal methods and verified computing produces rigorous estimates for trajectory computations with interval arithmetic, we still lack a definitive response to the request from the Radio Research Board for precise methods that determine the frequency of multiple time scale oscillators as functions of system parameters.

The second request from the Radio Research Board, mapping regions of parameter spaces with desired dynamical behaviors, is discussed in the next section of this paper on bifurcations. Despite the antipathy to graphical and numerical solutions expressed in the memorandum, identification of significant geometric structures via numerical solutions has become the predominant strategy for studying dynamical systems. Digital computers are capable of producing highly accurate simulations of trajectories that could hardly be imagined when the memorandum was written. Analyzing and interpreting the output of such computations is best done in a geometric setting. The \emph{phase portrait} of a vector field organizes all of its trajectories in its state space. Special types of trajectories, namely equilibrium points and periodic orbits\footnote{The \emph{equilibrium points} of the system $\dot{x} = f(x)$ are the solutions of $f(x) = 0$. A trajectory $x(t)$ is periodic with period $T$ if $x(T) = x(0)$ but $x(t) \ne x(0)$ for $0<t<T$}, are significant entities in phase portraits. Trajectories that tend to one of these special orbits in forward/backward time comprise its stable/unstable manifold. An equilibrium is \emph{hyperbolic} if the real parts of its Jacobian derivative eigenvalues are all non-zero. Periodic orbits are hyperbolic if none of the eigenvalues of a return map derivative have magnitude one. Stable and unstable manifolds of hyperbolic special orbits are indeed manifolds. These objects together comprise the phase portraits of planar vector fields, but Cartwright and Littlewood showed the presence of much more complicated dynamics in \eJG{FVDP}. 

Bifurcation analysis seeks to partition the parameter spaces of families of vector fields into equivalence classes with qualitatively similar dynamics. Bifurcations happen at parameters that are not in the interior of their equivalence classes. Ideally, the equivalence classes consist of topologically equivalent systems, but that viewpoint has aspects that fit the memorandum's precaution about solutions that are ``so complicated a form" that they fail be useful in practical settings.  A starting point that avoids these complications is to consider equivalence classes consisting of systems with the same numbers of equilibria and periodic orbits and the same dimensions of their stable and unstable manifolds. Bifurcations in this setting have explicit \emph{defining} equations (and inequalities) that locate submanifolds of the parameter space where specific types of bifurcations occur. See the next section for more details. 

Nonlinear dynamics is an interdisciplinary research area that has long addressed the ``deficiency of liaison'' between mathematics and other disciplines noted by the Radio Research Board memorandum. The boundaries between disciplines have become fuzzier over time, but the incisiveness of disciplinary research remains important.  The strategy adopted by Cartwright and Littlewood of formulating and studying mathematics problems based upon discoveries in other disciplines has repeatedly reaffirmed Wigner's widely quoted observation of ``the unreasonable effectiveness of mathematics in the natural sciences"~\cite{MR824292}. The synergy between mathematical studies of dynamical systems and experimental research in other disciplines has been remarkable. Observations of dynamical phenomena in empirical data that are not (yet) explained by existing theory have stimulated mathematical investigations yielding results of great generality. Conversely, mathematical discoveries about dynamical behaviors have prompted experiments that would not have been attempted otherwise. Studies of the ``transition to turbulence" in fluid systems~\cite{RT71,MR629208} are a notable example of this feedback between experiment and theory.  

New ``deficiencies of liaison" between disciplines continue to appear, even on issues that demand our attention. Two timely  examples today are 
\begin{itemize}
\item
Control of the electric grid~\cite{NAP21919}:   Alternating current electric grids operate as oscillators, but on a vastly larger scale than radios. Control of the grids to maintain constant voltage and frequency requires a continuing balance between power production and loads. Wind turbines and photovoltaic solar panels as renewable sources of electric power generation are intermittent, significantly complicating the control problems and raising new mathematical questions. 
\item
Climate ``tipping points"~\cite{tipping2023}: Recent increases in extreme weather and sea level rise caused by melting ice have raised concern about threats to humanity caused by rapid changes in climate systems. Studies of low dimensional multiple time scale dynamical systems have developed diagnostic tests to determine whether a system is approaching a tipping point. Melting permafrost is one climate system being analyzed with these tests. When permafrost melts, it releases methane into the atmosphere. Since methane is a potent greenhouse gas, its release causes further warming that could trigger tipping points called ``compost bombs"~\cite{MR2782155}. A second example of a potential climate tipping point is abrupt reorganization of Atlantic Ocean circulation accompanied by a much colder climate in northern Europe as happened in past ice ages.
\end{itemize}

\eJG{The Radio Research Board memorandum and much of the subsequent research by Cartwright and Littlewood on nonlinear oscillations concerns ``small" nonlinearities which can be studied with asymptotic methods that compute approximations in terms of a small parameter that represents the magnitude of the nonlinearity. However, the regime of ``relaxation oscillations" discussed in the 1945 paper of Cartwright and Littlewood~\cite{CL45} and the 1957 papers of Littlewood~\cite{Littlewood57,Littlewood57b} is quite different.  The remainder of this paper discusses the multiple time scale regime of relaxation oscillations.}

\section{Bifurcation Theory of Dynamical Systems}

Bifurcation theory characterizes how invariant sets change topologically in generic families of dynamical systems. In its simplest manifestation, bifurcation analysis counts the number of equilibrium points and periodic orbits and the dimensions of their stable and unstable manifolds. In the family of systems $\dot{x} = f(x,\lambda)$ with $x \in \R^n$ as state space variable and $\lambda \in \R^l$ as parameter, the equilibrium points form a smooth manifold of dimension $l$ near points where $Df = (D_x f,D_\lambda f)$ has full rank. Near points where $D_x f$ also has full rank, $x$ can be written as a function of $\lambda$ on this manifold. However if $D_x f$ is singular and suitable inequalities are satisfied by second derivatives of $f$, the manifold has a \emph{fold}: its projection to the parameter space has a boundary with two nearby equilibria on one side of the boundary and none on the other. The pair of \emph{defining} equations $f(x,\lambda) = 0, \, \det(D_x f)(x,\lambda) = 0$ locate fold bifurcations. If this pair of equations is regular with full rank $n+1$ and satisfies an additional inequality, projection of the set of nearby fold bifurcations to the parameter space is a manifold of dimension $(l-1)$ or \emph{codimension} $1$. When $x,\lambda \in \R$ are scalar variables, the needed inequalities are just $f''(x) \ne 0$ and $\frac{\partial f}{\partial \lambda} \ne 0$. The equation $\dot{x} = x^2 + \lambda$ is a simple example of the bifurcation, a \emph{normal form} that embodies the qualitative changes that occur at this bifurcation. Fold bifurcations for periodic orbits of $f$ are characterized by analogous defining equations for fixed points of a return map of $f$.  
 
 Defining equations can also characterize bifurcations of higher codimension. Consider the family of equations $\dot{x} = x^3 + \lambda_1 x + \lambda_2$ with one dimensional state space and two dimensional parameter space.  The  equilibrium points of this system are given by the smooth surface $x^3 + \lambda_1 x + \lambda_2 = 0$ in the product space. The projection of this surface onto the $(\lambda_1,\lambda_2)$ parameter space has a \emph{cusp}: it is singular   when $3x^2 + \lambda_1 = 0$. Non-zero points on this curve are fold bifurcations, but the origin fails to satisfy the inequality $f''(x) \ne 0$ that is part of the characterization of a fold bifurcation. Instead, the three equations $f=0,\, f' = 0, \, f'' = 0$ (with additional inequalities involving the parameters) locate \emph{cusp} bifurcations of codimension 2 in the parameter space of a system with one dimensional state space.
 
Classifying bifurcations with defining equations and \emph{nondegeneracy} inequalities produces a dictionary of different types of bifurcations that appear on submanifolds in the parameter spaces of generic systems~\cite{MR1736704,MR4592585}. Together, the bifurcations located by this process form a \emph{stratified set} in which codimension $(k+1)$ bifurcations occur on boundaries of codimension $k$ bifurcations. In some cases, codimension $(k+1)$ bifurcations are in the boundary of multiple strata of codimension $k$ bifurcations, including \emph{global} bifurcations that are not captured by just counting equilibrium points and periodic orbits. Analysis of families that are transverse to the bifurcation manifolds (called \emph{unfoldings} or \emph{normal forms}) reveals these global bifurcations. The stratified sets of bifurcations located by these procedures yield a partial dictionary of geometric structures that may be encountered in mapping the parameter spaces of generic systems.

We end this section with a brief description of bifurcations that are not found by the process described above. 
The family of planar diffeomorphisms $h(x,y) = (1-ax^2 + y,bx)$, studied first by  H\'enon~\cite{Henon76},  has a horseshoe when $b>0$ is small and $a > 2+b$ while when $a >>b>0$, and $a$ is small, all bounded trajectories tend to just one of two fixed points. Thus an infinite number of bifurcations occur as $a$ varies. When $b=0$, the H\'enon family reduces to the family of non-invertible quadratic maps of the interval. The dynamics of the quadratic family of one dimensional maps have been characterized in remarkable detail, including 
\begin{itemize}
\item
the order in which the bifurcations of periodic orbits occur~\cite{MR970571,MR438399}, and 
\item
the existence of a nowhere dense parameter set of positive Lebesgue measure for which the maps have \emph{chaotic attractors} with absolutely continuous invariant measures~\cite{MR630331,MR799250}. 
\end{itemize}
The existence of chaotic attractors has also been established for the H\'enon map with $b>0$~\cite{BC91,MV93,WY01}.They are not structurally stable, but their chaotic attractors are a better match to those encountered in systems like FVDP  than those in structurally stable systems. 

\section{Geometric Singular Perturbation Theory}

\emph{Geometric singular perturbation theory}~\cite{Fenichel77,Jones95,MR3309627} (acronym GSPT) studies dynamical systems with multiple time scales. Many dynamical phenomena are greatly distorted in systems with multiple time scales, giving rise to distinctive aspects of GSPT.  The subject is complicated geometrically, analytically and computationally. It remains incomplete and has not been widely disseminated even within the dynamical systems community.

Most of GSPT created thus far applies to  \emph{slow-fast} vector fields that take the form
\begin{align}
\eps \dot u & = f(u,v) \\
\dot v & = g(u,v)
\label{sf}
\end{align}
where $u \in R^k$ are the \emph{fast} variables, $v \in R^m$ are the \emph{slow} variables and $\eps > 0$ is a small parameter representing the ratio of time scales. The \emph{singular limit} $\eps = 0$, is a \emph{differential algebraic equation} (acronym DAE) whose analysis is a starting point for analysis of the \emph{full} system with $\eps > 0$.  
The set $f(u,v)=0$  is the \emph{critical manifold} $C$ of \eqref{sf}. Near points of $C$ where $D_vf$ has full rank $m$, 
$C$ can be written as the graph of a function $u = h(v)$ and the DAE reduces to the ODE $\dot{v} = g(h(v),v)$, called the \emph{slow flow}. At points of $C$ where $D_vf$ is singular, the existence and uniqueness theorem for ordinary differential equations fails. This is evident at the origin for the DAE $v+u^2 = 0;\,\dot{v} = 1$ where $v= -u^2$ is inconsistent with $\dot{v} = 1$. A change of time scale in equation~\eqref{sf} produces the \emph{fast subsystem} or \emph{layer equation} $\dot{u} = f(u,v)$ in which $v$ acts as a parameter and the critical manifold is comprised of equilibrium points. The critical manifold is \emph{normally hyperbolic} at points $(u_0,v_0)$ where $D_vf$ has full rank $m$ and the equilibrium points of the layer equation $\dot v = g(u_0,v)$ are hyperbolic, i.e. $D_vg$ has eigenvalues with non-zero real parts. 

The foundational theorem of GSPT~\cite{Fenichel79,Jones95} establishes the existence of locally invariant \emph{slow manifolds} of a slow-fast system that approximate the slow flows of normally hyperbolic critical manifolds. This theorem is complicated in several respects. First, the slow manifolds are parameterized by the time scale parameter $\eps$ and only have meaning in the asymptotic setting of $\eps \to 0$. Second, slow manifolds are ``locally invariant," in that they contain only finite time segments of trajectories. Trajectories leave slow manifolds through their boundaries, but the location of these boundaries is a matter of choice. Third, invariant manifolds that are ``exponentially close" (i.e, of the order of $\exp(-c/\eps),\, c> 0$) to a slow manifold are also slow manifolds. Fourth, iterative numerical algorithms to compute normally hyperbolic slow manifolds have been proposed, but need to specify boundary conditions carefully in order to converge to a specific slow manifold.

\begin{figure}[t!]
\centering
\includegraphics[width=4in]{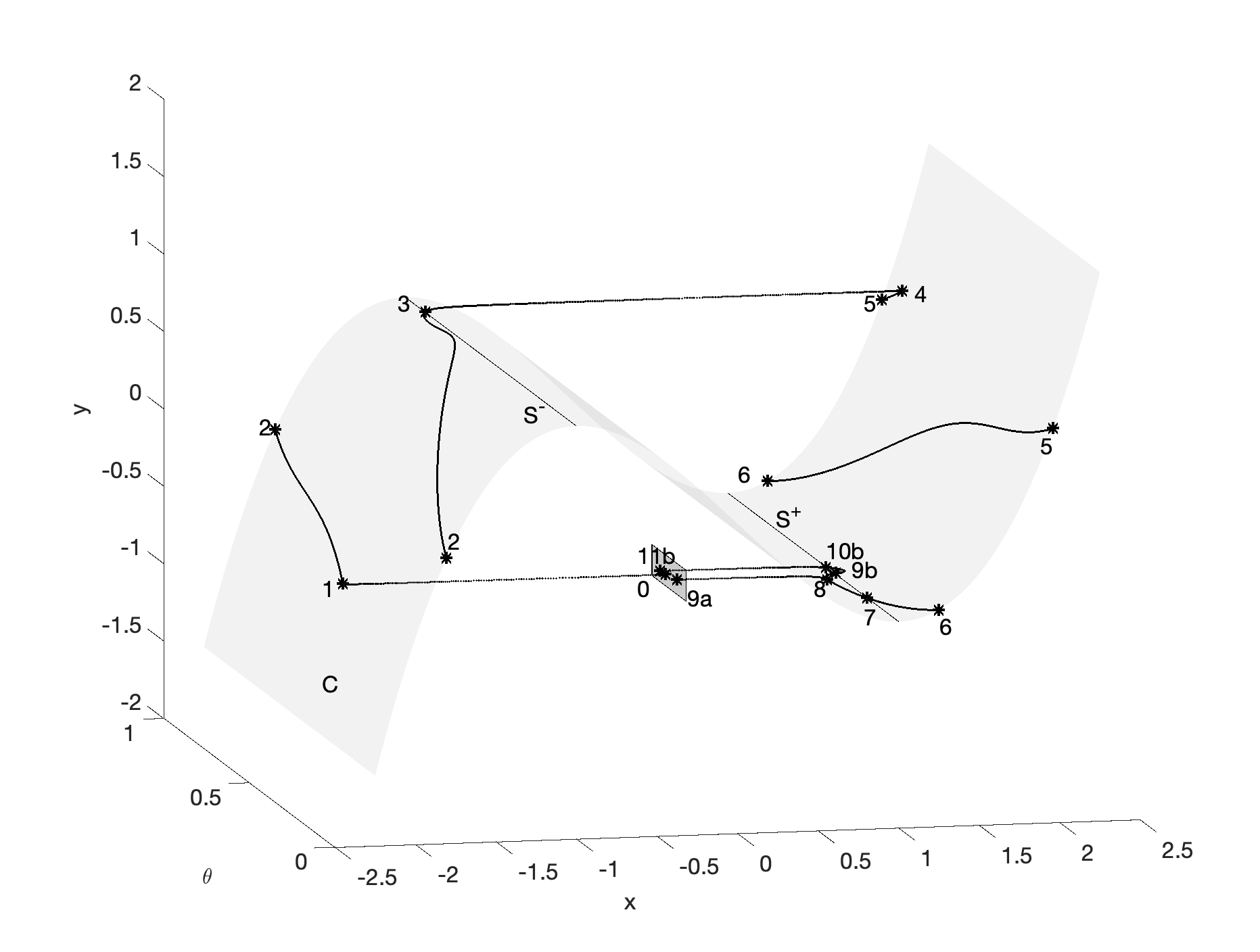}

\caption{\label{fvdp_pp} The critical  manifold of FVDP and returns of two trajectories to a cross-section defined by $x=0$: Parameters in FVDP are $(a,\omega,\eps) = (1.1,1.505,0.001)$. The state space variable $\theta$ lies in $\R/\Z$, so points with $\theta=0$ and $\theta=1$ are identified. The critical manifold $C$ is displayed in light gray and the cross-section in a slightly darker shade of gray. Significant points on the trajectories are numbered and marked with asterisk symbols. The two trajectories have initial points at $(0,-0.6752,0.41732694)$ and $(0,-0.6752,0.41732695)$ near the point labelled $0$ and are visually indistinguishable till they reach the point labelled $8$. To get there, they jump (flow in the fast direction almost parallel to the $x-$axis) until they cross $C$ at the point labelled $1$. There, they turn and flow along the slow manifold until they reach $x = -1$ at the point labelled $3$. The points labelled $2$ mark the location where $\theta$ reaches the value $1$ and ``resets" to $0$. The trajectory segment from $3$ to $4$ is a fast jump terminating at the slow manifold with $x>1$. The flow along the slow manifold with $x>1$ reaches $\theta = 1$ and resets to $\theta = 0$ twice at the points labelled $5$ and $6$. At $7$ the trajectories cross $x=1$ where they continue together as canards close to the repelling sheet of $C$. Near the point labelled $8$, the trajectories abruptly jump in opposite fast directions. One trajectory, the ``slice," jumps toward the cross-section $x=0$, reaching it at the point labelled $9a$. The second trajectory, the ``dip," jumps back to the slow manifold with $x > 1$, crossing $x=1$ at point $9b$. This trajectory makes a third crossing of $x=1$ at $10b$ and then jumps again, reaching the cross-section $x=0$ at the point labelled $11b$.
}
\end{figure}

FVDP is a slow-fast system with slow variables $(y,\theta)$ and fast variable $x$. Its critical manifold is the cubic surface $C$ defined by $y + x - \frac{x^3}{3} = 0$. 
$C$ is normally hyperbolic except on the \emph{fold curves} $S^{\pm} = \{(x,y,\theta)| x = \pm 1, y = \mp \frac{2}{3})\}$ where $C$ has a tangent vector parallel to the fast $x$ direction. Near points of $C - S^{\pm}$, $C$ can be written as the graph of a function $x = h(y,\theta)$ and the DAE reduces to the slow flow, defined by $\frac{d}{dt}(y,\theta) = (-h(y,\theta) + a \sin(2\pi \theta), \omega)$. 
Figure \ref{fvdp_pp} shows the critical manifold of FVDP and two trajectories with initial conditions close to one another. 

The slow flow of FVDP can be extended to the fold curves of $C$ by employing $(x,\theta)$ coordinates and rescaling time by a factor $(x^2-1)$, yielding the \emph{desingularized slow flow}
\begin{equation}
\begin{array}{rcl}
  \dot{x} & = & -x  + a \sin(2\pi\theta)\\
  \dot{\theta} & = & \omega(x^2-1)
 \end{array}
\label{fvdpvd}
\end{equation}
Note that 
\begin{itemize}
\item
the time rescaling reverses the orientation of the slow flow in the strip $-1 < x < 1$, 
\item
the desingularized slow flow has equilibrium points (called \emph{folded equilibria}) at $(x,\theta) = (\pm 1,\sin^{-1}(\frac{\pm 1}{2 \pi a}))$ when $|a| > 1$, and 
\item
the desingularized slow flow is orthogonal to the fold curves on $C$. 
\end{itemize}
Phase portraits of the desingularized slow flow can be computed easily with standard numerical integration algorithms. 

In the unforced van der Pol vector field
\begin{equation}
\begin{array}{rcl}
   \eps \dot{x} & = & y + x - \frac{x^3}{3}  \\
   \dot{y} & = & - x 
 \end{array}
\label{vdp}
\end{equation}
on the $(x,y)$ plane, singular limits of trajectories concatenate segments of the critical curve $ y = \frac{x^3}{3} - x$ with horizontal segments from $(-1,\frac{2}{3})$ to $(2,\frac{2}{3})$ and from $(1,-\frac{2}{3})$ to $(-2,-\frac{2}{3})$. The  jumps along the horizontal segments are instantaneous on the slow time scale of the system. Limiting trajectories in the full forced van der Pol system make similar jumps when they cross the fold curves, as evident in Figure~\ref{fvdp_pp}. However, more complex dynamics occur near the folded equilibria. 

\begin{figure}[t!]
\centering
\includegraphics[width=4in]{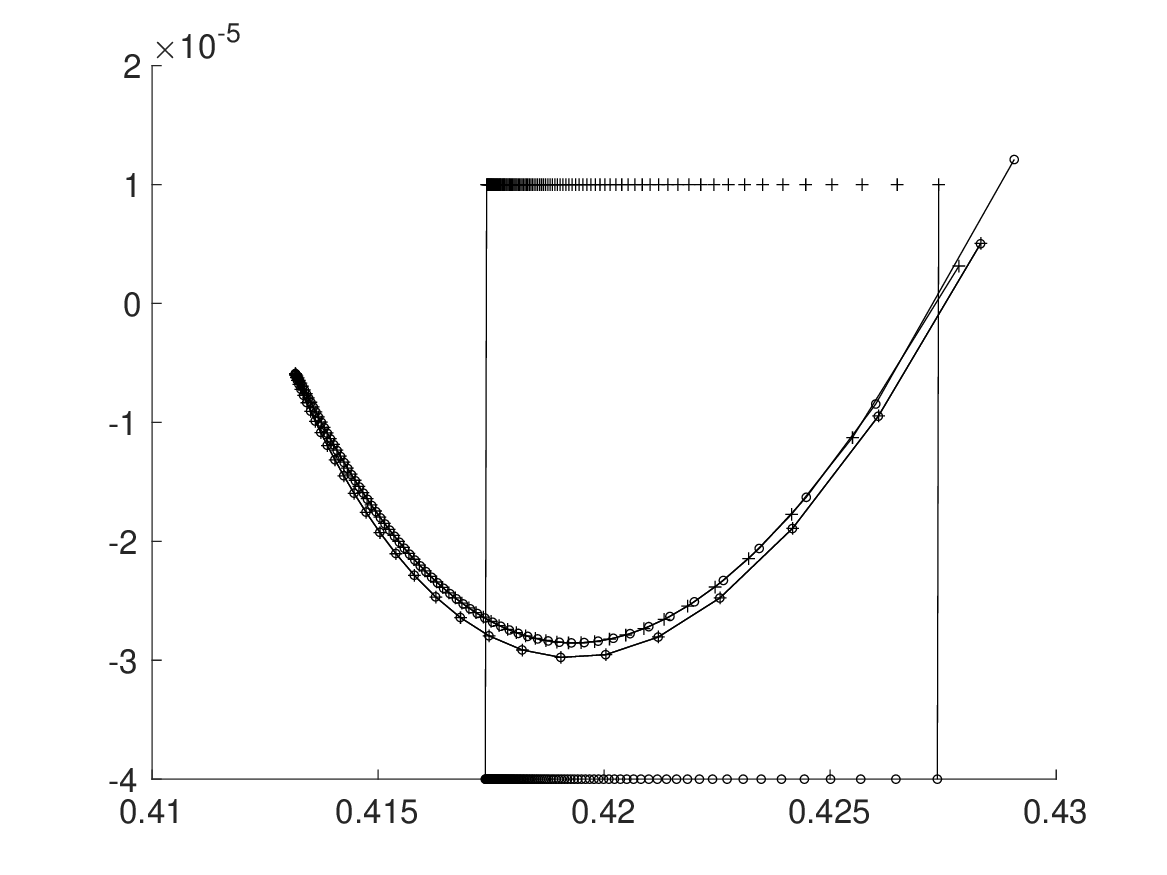}
\caption{ \label{fvdp_horse}
The image of a quadrilateral under a return map of FVDP. The system parameters are the same as those in Figure \ref{fvdp_pp}. The quadrilateral lies in the plane $x=0$ and has vertices $(0.4174, -0.67509216),(0.4274,-0.67657616),(0.41737,-0.675137708),( 0.42737,-0.676621708)$ in the $(\theta,y)$ plane close to the initial points of the two trajectories displayed in Figure 1. In this figure, the horizontal coordinate is $\theta$ and the vertical coordinate is chosen so that it is contracted by the return map. Images of points on the top and bottom of the quadrilateral under the return map to the plane $x=0$ with $x$ decreasing are plotted with symbols $+$ and $\circ$. The return map images of both top and bottom boundaries of the quadrilateral are curves that stretch across the quadrilateral from right to left, fold to the left of the quadrilateral, and then stretch back across the quadrilateral from left to right, thereby giving evidence that the quadrilateral creates a horseshoe of the return map. Note the disparity in the ranges of the horizontal and vertical coordinates: the image of the quadrilateral and its image are much wider in the $\theta$ direction than the transverse direction.
}

\end{figure}

The folded equilibria of the desingularized slow flow can be classified in terms of the eigenvalues of their Jacobians. Folded saddles with one positive and one negative eigenvalue play an especially important role in shaping the dynamics of FVDP. In the desingularized slow flow, the stable manifold of a folded saddle consists of two trajectories, one on the attracting sheet of the critical manifold and one on the repelling sheet. Some trajectories of the full system flow ``through" the folded saddle and follow its stable manifold onto the repelling sheet before abruptly jumping along the fast direction to one of the two attracting slow manifolds. The trajectory segments that follow the repelling sheet of $C$  for a distance $O(1)$ are called \emph{canards}: they create extreme divergence of nearby trajectories. The trajectories displayed in Figure~\ref{fvdp_pp} have canards after the trajectories pass near the $x=1$ fold curve of the critical manifold to reach the repelling slow manifold.

Canards were discovered and named by the Strasbourg ``duck hunters" \cite{Diener} who investigated slow-fast systems in terms of non-standard analysis.  One member of this group, \'Eric Beno\^it, analyzed folded saddles~\cite{Benoit83}. He proved that extensions of the slow manifolds near the attracting and repelling sheets of the critical manifold intersect along a ``maximal" canard trajectory at an angle that is $O(\eps)$. Nearby trajectories to the attracting slow manifold follow the maximal canard until they turn abruptly and jump away from the repelling slow manifold. The two trajectories displayed in Figure~\ref{fvdp_pp} jump in opposite directions: one jumps back to the slow manifold with $x > 1$ from which it came (a \emph{dip} in the terminology of Cartwright and Littlewood), and one jumps to the slow manifold with $x < -1$ (a \emph{slice} in the terminology of Cartwright and Littlewood.) The divergence of trajectories as they jump gives rise to horseshoes in the dynamics of FVDP. Figure~\ref{fvdp_horse} shows a small quadrilateral in the cross-section $x=0$ and its return as a ``horseshoe" to this cross-section.  Along the direction of the slow manifold, the image is stretched and folded. However, the folding is hardly apparent in the figure because the return is very thin in the cross-section.

\section{Numerical Methods and Results}

At a 1997-98 program in dynamical systems at the Institute for Mathematics and its Applications in Minneapolis, I began a collaboration with Kathleen Hoffman and Warren Weckesser that included new investigations of FVDP. We began by developing numerical methods to compute the slow flow on the critical manifold with high precision, together with its jumps at fold curves~\cite{GHW03}. One outcome of this study was a quantitative characterization of bifurcations in a \emph{hybrid} dynamical system that represents the slow flow, including trajectories that follow a canard before executing jumps away from the canard. We then extended this analysis to the full system with $\eps >0 $ small, paying special attention to canards and invariant slow manifolds~\cite{GHW03b}. Our work discovered several types of bifurcations that had not been described previously as well as parameter regimes where the system seemed to be structurally stable. We made essential use of boundary value solvers implemented in the computer program AUTO~\cite{AUTO}.\footnote{The matlab package MATCONT~\cite{Matcont} has similar capabilities as well as methods for locating and continuing codimension one and two bifurcations in multidimensional parameter spaces.}

Our conclusions about the dynamics of FVDP were not rigorous since they lacked bounds on errors in the computed data. This is a general issue for computational studies of dynamical systems that has been addressed with algorithms employing interval arithmetic. One step algorithms for solving initial value problems of ordinary differential equations march along a trajectory, iteratively computing approximations at successive time steps. These computations incur truncation errors due to the approximation of the solution and rounding errors from the use of floating point arithmetic. Verified computing establishes rigorous bounds for these errors, but the bounds inherently grow for trajectories that have a positive Lyapunov exponent -- as is the case for the FVDP horseshoes. To circumvent this difficulty, computer proofs that confirm the existence of horseshoes construct and verify the existence of a \emph{hyperbolic structure} in a domain containing the horseshoe and then use the dynamical systems concept of \emph{shadowing}~\cite{MR2423393} which establishes long time properties of trajectories from short time approximations. The shadowing property states that close to a curve whose short time segments have small errors, there is a unique trajectory. The distance from a numerically computed trajectory to the ``true" trajectory with the same initial condition becomes large, but there is another initial condition whose trajectory remains close to the computed one. Haiduc~\cite{Haiduc09} used verified estimates of short time trajectory segments of the slow flow and the existence of a hyperbolic splitting in the full system for small enough values of the time scale parameter to produce a rigorous proof for the existence of hyperbolic invariant sets for FVDP. He further showed that the rest of the nonwandering set consists of a single repelling periodic orbit and two stable periodic orbits. These results constitute a satisfying end to the ``preliminary results" of Cartwright and Littlewood in their 1945 paper by confirming that FVDP has parameter regions with chaotic invariant sets in which it is structurally stable. 

What can be said about the dynamics of FVDP in other parameter regions? In the $(a,\omega)$ parameter plane of FVDP with $a>2$, there are strips in which the system has stable periodic orbits of odd period. The strips with periods $(2n\pm1)$ overlap and chaotic basic sets are found in these overlaps. Boundaries of the strips were computed with asymptotic methods by Flaherty and Hoppensteadt~\cite{MR499449}. There are also parameter regions in the strips with a single stable periodic orbit where the only basic sets are a few periodic orbits. Thus, bifurcations like those found in the H\'enon map must occur along parameter curves that connect these two types of regions. This suggests that there are chaotic attractors at some parameter values along these paths. While that has not yet been proved, Guckenheimer, Wechselberger and Young~\cite{GWY06} proved that such attractors do occur in a forced relaxation oscillator created by modifying the  FVDP equations.

\section{More Applications}
The electronic devices modeled with FVDP by van der Pol have been supplanted by solid state electronics, diminishing the incentive to pursue further detailed study of this particular model. However, many other applications have used dynamical systems models and raised additional mathematical questions. Here are brief descriptions of four of these applications:
\begin{itemize}
\item
Semiconductor lasers with optical injection have been modeled by rate equations that are three dimensional vector fields, the same size as FVDP. These devices exhibit many different dynamical states including ones that are chaotic. Detailed comparisons between experimental results varying forcing frequency and amplitude with bifurcation analysis of the models show remarkable agreement~\cite{wieczorek2005}. Indeed, previously unobserved states first seen in a model were subsequently found in new experiments. This system is a role model for the use of bifurcation analysis to make quantitative predictions about complicated dynamical behaviors in physical systems. It has also stimulated research on differential delay equations since delays are important in many laser systems.
\item
The Belousov-Zhabotinsky reaction in stirred tank reactors consists of a mixture of chemicals that is capable of chaotic oscillations.\footnote{Because there is a flow of reactants through the reactor, the system is not constrained by the second law of thermodynamics to tend to an equilibrium state.} Whether these chaotic oscillations are the result of noise or whether they are an inherent feature of idealized reactors in noise free environments was a controversial issue during the 1980's.  The stiffness\footnote{\emph{Stiffness} is a commonly used term without a precise meaning in reference to simulation of systems of ordinary differential equations that require \eJG{a very large number of} time steps to obtain meaningful results \eJG{with \emph{explicit} integration algorithms}. Integration of systems that have a wide range of time scales for times of the order of the slowest time scales are stiff. In particular, kinetic equations for chemical reactors often have reaction rates that differ by ten orders of magnitude or more.} of the more detailed models precluded accurate numerical simulations. Several distinct types of complicated dynamics have been studied experimentally~\cite{zhang93}, including \emph{mixed mode oscillations}~\cite{MR2916308} in which concentrations switch repeatedly between small and large amplitudes. Guckenheimer and Scheper~\cite{MR2788920} studied two types of mixed mode oscillations in a nine dimensional differential equation model based on mass action kinetics and gave a rough comparison between the simulation results and a theoretical classification of dynamical mechanisms that produce mixed mode oscillations. This is an unfinished story in all of its aspects. Further work to produce quantitative fits of model dynamics with experimental data from the Belousov-Zhatbotinsky might well raise new mathematical questions.
\item
Conductance based models of neuronal systems are electrical circuits that represent membranes as capacitors containing embedded channels that selectively allow ionic currents of one or more species to flow across the membrane when open~\cite{10.5555/94605}. The current through each family of channels is modeled as a nonlinear resistor whose conductance typically depends upon voltage across the channel and the binding of neuromodulators. The force driving these currents comes from ionic concentration differences across the membrane maintained by slower processes. Computational neuroscience creates dynamical systems models whose state space variables are membrane potentials, gating variables of different families of voltage dependent ionic conductances and sometimes relevant ion concentrations (especially $Ca^{++}$) as models of electrical activity in the nervous system. The 1952 Hodgkin-Huxley model of squid giant axon~\cite{Hodgkin1952} is a system of four differential equations that is the archetype of conductance based models. It models action potentials, rapidly activating spikes in membrane potential whose propagation is a signaling mechanism of nervous systems: action potentials in presynaptic neurons trigger synaptic currents in postsynaptic neurons. More complex temporal patterns of action potentials that undergo bifurcations to different behavioral states are observed in different systems. One system among many that have been studied extensively is the central pattern generator for respiration~\cite{doi:10.1152/jn.90958.2008}. Since gating of membrane currents occurs over a broad range of time scales and varied rhythmic behaviors are observed, GSPT has had a big impact. Empirical time series measurements have raised new mathematical questions while studies of conductance based models have generated hypotheses to be tested experimentally. The interactions between mathematics and biology in computational neuroscience have enriched both subjects greatly.
\item
There has been a resurgence in the study of pattern formation within infinite dimensional dynamical systems defined by partial differential equations. The \emph{FitzHugh-Nagumo equation} is a simplified model of action potential propagation along a nerve axon as a nonlinear reaction diffusion equation. Traveling wave solutions to the equation are described by trajectories of an ordinary differential equation reduction of the model in which the wave speed enters as a parameter. In some parameter regimes, the reduced model is a slow-fast vector field  with one slow variable and two fast variables. It inspired the discovery of a seminal result in GSPT, the \emph{Exchange Lemma}, by Jones and Kopell~\cite{MR1268351}. Tangencies of invariant manifolds result in enigmatic bifurcation behavior of this system whose analysis required further extensions of GSPT and numerical methods for computing invariant slow manifolds of saddle type~\cite{MR2595892,MR3755655}.
\end{itemize}

\section{A Concluding Remark}

The 1945 Cartwright-Littlewood ``preliminary survey" of their research on the forced van der Pol equation had far more impact than the detailed proofs of their results published by Littlewood in 1957. Rigorous proofs as the heart of mathematics is a deeply entrenched community value, but the 1945 Cartwright-Littlewood paper exemplifies that fostering -- and evaluating -- creative ``preliminary" research is also important to the discipline. The memorandum of the Radio Research Board that sparked their collaboration and the subsequent development of ``chaos theory" also highlight the complementarity of quantitative studies of specific dynamical models for applications and qualitative research that connects seemingly unrelated phenomena. Cartwright and Littlewood left a lasting legacy on both of these pursuits.

\section*{Acknowledgment}
I give special thanks to Caroline Series, June Barrow-Green and the referees who helped improve the historical accuracy of this paper

\bibliographystyle{plain}
\bibliography{CL}

\end{document}